\documentclass[a4paper,11pt]{article}

\usepackage{latexsym}
\newtheorem{thm}{Theorem}[section]

\newenvironment{pf}{\noindent\textbf{Proof\ }}{\hfill$\Box$\smallskip}
\newtheorem{cor}[thm]{Corollary}
\newtheorem{prop}[thm]{Proposition}
\newtheorem{conj}[thm]{Conjecture}
\newtheorem{rem}[thm]{Remark}
\newcommand{\ov}{\overline}
\newcommand{\om}{\omega}
\newcommand{\F}{\mathbf{F}}
\newcommand{\mx}[4]{\left(
\begin{array}{cc}
{#1}&{#2}\\
{#3}&{#4}
\end{array}
\right)}
\newcommand{\vc}[2]{\left(
\begin{array}{c}
{#1}\\
{#2}
\end{array}
\right)}
\newcommand{\GL}{\mathop{\mathrm{GL}}}

\newcommand{\ce}{\mathcal{E}}
\newcommand{\Z}{\mathbf{Z}}
\newcommand{\N}{\mathbf{N}}
\newcommand{\ti}{\tilde}
\newcommand{\al}{\alpha}
\newcommand{\be}{\beta}
\newcommand{\ga}{\gamma}

\newcommand{\tht}{\theta}

\title{Symmetric Pascal matrices modulo p}
\author{Roland Bacher\footnote{Support 
from the Swiss National Science Foundation is
gratefully acknowledged.} and Robin Chapman}
\date{30 January 2003}
\begin{document}
\maketitle

\section{Introduction}

This paper presents results and conjectures concerning 
symmetric matrices associated to Pascal's triangle.
We first give a formula for the determinant over $\Z$ of the reduction 
modulo $2$ with values in $\{0,1\}$ for such a matrix. We then study
the reduction modulo a prime $p$ of the characteristic polynomials of these 
matrices. Our main results imply a formula for the prime $p=2$ and
a conjectural formula for $p=3$.

Consider the symmetric matrix $P(n)$ with coefficients
$$p_{i,j}={i+j\choose i},\ 0\leq i,j<n\ .$$
We call $P(n)$ the \emph{symmetric Pascal matrix} of order~$n$.
The entries of $P(n)$ satisfy the recurrence
$$p_{i,j}=p_{i-1,j}+p_{i,j-1}.$$
In \cite{B} the first author studied the determinant of the general
matrix with entries satisfying this recurrence.

An easy computation yields $P(\infty)=T\ T^t$ where $T$ is the infinite
unipotent lower triangular matrix
$$T=\left(\begin{array}{cccccccccc}
1\cr
1&1\cr
1&2&1\cr
1&3&3&1\cr
\vdots&&&&\ddots\end{array}\right)=\hbox{exp}\left(\begin{array}{cccccccccc}
0\cr
1&0\cr
0&2&0\cr
&0&3&0\cr
&&&&\ddots\end{array}\right)$$
with coefficients $t_{i,j}={i\choose j}$. This shows that $\det (P(n))=1$
and that $P(n)$ is positive definite
for all $n\in {\bf N}$. Hence all zeroes
of the characteristic
polynomial $\chi_n(t)=\det(tI(n)-P(n))$
(where $I(n)$
denotes the identity matrix of size $n$) of $P(n)$ 
are positive reals. The inverse $P(n)^{-1}$ of $P(n)$ 
is given by 
$$P(n)^{-1}=\left(T(n)^t\right)^{-1}T(n)^{-1}$$
and $T(n)^{-1}$ has coefficients $(-1)^{i+j}{i\choose j},\ 0\leq i,j<n$.
Hence $T(n)$ and $T(n)^{-1}$ are conjugate, and thus also
$P(n)$ and $P(n)^{-1}$ are conjugate. The 
characteristic polynomial $\chi_n(t)$ therefore satisfies
$\chi_n(t)=(-t)^n\chi(1/t)$
and $1$ is always an eigenvalue of $P(2n+1)$, cf. \cite{L}.
The polynomials $\chi_n(t)$, especially their behaviour modulo
primes, will be our main object of study. For convenience, we
write $I$ for $I(n)$ whenever the size of the identity
matrix is unambiguous.

Define ${\overline P}(n)_2$ as the reduction modulo $2$ of $P(n)$
with values in
$\{0,1\}$ by setting
$${\overline p}_{i,j}=\left({i+j\choose i}\pmod 2\right)\in\{0,1\}\ .$$

The Thue-Morse sequence $s_n=\sum\nu_i\pmod 2$ counts the parity of 
all non-zero digits of a binary integer $n=\sum \nu_i2^i$. It can also 
be defined recursively by $s_0=0$, $s_{2k}=s_k$ and $s_{2k+1}=1-s_k$
(cf. for instance \cite{AS}).

\begin{thm} \label{detmod2}
The determinant (over $\bf Z$) of ${\overline P}(n)_2$ is
given by
$$\det ({\overline P}(n)_2)=\prod_{k=0}^{n-1} (-1)^{s_k}\ .$$
\end{thm}

A similar result holds for the reduction modulo~$3$ of $P(n)$ 
with values in $\{-1,0,1\}$.

In the sequel, we will be interested in the characteristic polynomial
$\det(tI-P(n))\pmod p$ for $p$ a prime number. The next result yields
a formula for $n=p^l$ and is of crucial importance in the sequel.

\begin{prop}\label{order3forprimep} 
Given a power $q=p^l$ of a prime $p$, the matrix $P(q)$
has order $3$ over $\F_p$. Its characteristic polynomial
$\chi_q(t)=\det (tI(q)-P(q))$ satisfies
$$\chi_q(t)\equiv(t^2+t+1)^{\frac{q-\epsilon(q)}{3}}
(t-1)^{\frac{q+2\epsilon(q)}{3}}\pmod p$$
where $\epsilon(q)\in \{-1,0,1\}$ satisfies $\epsilon(q)\equiv q\pmod 3$.
\end{prop}
In particular, $P(q)$ can be diagonalized over $\F_{p^2}$ except when $p=3$.
For instance, $P(3)$ has a unique Jordan block over $\F_3$.

This proposition (except for the diagonalization part) admits the following 
generalization:

\begin{thm}\label{ppowerformula} 
When $q=p^l$ is a power of a prime~$p$ and
$0\leq k\leq q/2$ then
$$\chi_{q-k}(t)\equiv
(t^2+t+1)^{(q-\epsilon(q))/3-k}(t-1)^{(q+2\epsilon(q))/3-k}
\det(t^2 I+P(k))\pmod p$$
where $\epsilon(q)\in \{-1,0,1\}$ satisfies $\epsilon(q)\equiv q\pmod 3$.
\end{thm}

Theorem \ref{ppowerformula} completely determines the
reduction modulo $2$ of $\chi_n(t)$ as follows:
Define a sequence $\gamma(0)=0,
\gamma(1),\dots$ recursively by
$$\gamma(2^l-k)=\frac{2^l+2(-1)^l}{3}-k+2\gamma(k),\ 0\leq k\leq
2^{l-1}\ .$$

\begin{thm} \label{mod2} 
For all $n\in \N$
$$\chi_n(t)\equiv (t+1)^{\gamma(n)}
(t^2+t+1)^{\gamma_2(n)}\pmod 2$$
where $\gamma_2(n)=\frac12(n-\gamma(n))$.
\end{thm}

It follows immediately that the matrix $I-P(n)^3$
is nilpotent over $\F_2$ for all $n\in {\bf N}$.

The first terms $\gamma(1),\dots,\gamma(32)$ and $\gamma_2(1),
\dots,\gamma_2(32)$ are given by
$$\begin{array}{| c | cccccccccccccccc |}
\hline
n&1&2&3&4&5&6&7&8&9&10&11&12&13&14&15&16\cr
\gamma(n)&1&0&3&2&5&0&3&2&5&0&11&6&9&4&7&6\cr
\gamma_2(n)&0&1&0&1&0&3&2&3&2&5&0&3&2&5&4&5\cr
\hline
n&17&18&19&20&21&22&23&24&25&26&27&28&29&30&31&32\cr
\gamma(n)&9&4&15&10&21&0&11&6&9&4&15&10&13&8&11&10\cr
\gamma_2(n)&4&7&2&5&0&11&6&9&8&11&6&9&8&11&10&11\cr
\hline\end{array}$$

The sequence $\gamma(0),\gamma(1),\dots$ has many interesting
arithmetic features. In order to describe them, let us introduce
the number $b(n)$ of ``blocks''
of adjacent ones in the binary representation of a positive integer $n$.
For instance
$667=(1010011011)_2$ and so $b(667)=4$. 
Notice that $b(2n)=b(n)$ and
$b(2n+1)=b(n)+1-\left(n\pmod 2\right)$
(with $n\pmod 2\in \{0,1\}$). This, together with $b(0)=0$, defines
the sequence $b(n)$ recursively.

\begin{thm} \label{gammaformula}

\ \ (i) We have
$$\gamma(2^l+k)=\frac{2^l+2(-1)^l}{3}-k+4\gamma(k)$$
for all $0\leq k\leq 2^{l-1}$.

\ \ (ii) We have for all $n\in\N$ and $2^{l-2}\leq k\leq 2^{l-1}$
$$\gamma(2^l-k)=\gamma(k)+2\gamma(2^{l-1}-k)\ .$$

\ \ (iii) We have 
$$\gamma(2^l+k)=1+\gamma(2^l+k-1)+2\gamma(2^l-k)-2
\gamma(2^l+1-k)$$
for $1\le k\le 2^l$.

\ \ (iv) We have
$$\begin{array} {lcl}
\gamma(2n)&=&n-\gamma(n)\ ,\cr
\gamma(2n-1)&=&\gamma(2n)+(4^{b(2n-1)}-1)/3
=n-\gamma(n)+(4^{b(2n-1)}-1)/3\ ,\cr
\gamma(2n+1)&=&\gamma(2n)+(2^{1+2b(n)}+1)/3
=n-\gamma(n)+(2^{1+2b(n)}+1)/3\ .\end{array}$$
\end{thm}


Part (iv) of this Theorem gives an alternative recursive definition
of the sequence $(\ga(n))$.

Theorem \ref{ppowerformula} seems to have many generalizations.
A first one is given by the following:

\begin{conj}\label{firstgeneralization}
For each integer $k\ge0$ there exists a monic polynomial
$c_k(t)\in\Z[t]$ of degree $4k$ such that 
$c_k(t)=t^{4k}c_k(t^{-1})$ with the following property:
if $q$ is a power of a prime $p$, and $0\le k\le q/2$ then
$$\chi_{q+k}(t)\equiv(t^2+t+1)^{(q-\epsilon(q))/3-k}
(t-1)^{(q+2\epsilon(q))/3-k}c_k(t)\pmod{p}$$
where $\epsilon(q)\in \{-1,0,1\}$ satisfies $\epsilon(q)\equiv q\pmod 3$.
\end{conj}

The first few of these conjectural polynomials $c_k(t)$ are
\begin{eqnarray*}
c_0(t)&=&1,\\
c_1(t)&=&t^4-2t^3-2t+1,\\
c_2(t)&=&t^8-6t^7+4t^6-4t^5+15t^4-4t^3+4t^2-6t+1,\\
c_3(t)&=&(t^4-2t^3-2t+1)(t^8-16t^7+4t^6-4t^5+40t^4-4t^3+4t^2-16t+1),\\
c_4(t)&=&t^{16}-58t^{15}+288t^{14}-240t^{13}+393t^{12}-1440t^{11}+836t^{10}
-902t^9\\
&&{}+2376t^8-902t^7+\cdots-58t+1,\\
c_5(t)&=&c_1(t)(t^{16}-196t^{15}+2112t^{14}-792t^{13}+1290t^{12}-10560t^{11}\\
&&{}+2768t^{10}-2972t^9+17424^8-2972t^7+\cdots-196t+1).\\
\end{eqnarray*}



For $p=2$, it follows from Theorem \ref{mod2} and assertion (ii) in
Theorem \ref{gammaformula} that if $c_k(t)$ exists then
$$c_k(t)\equiv \left(\det(tI+P(k))\right)^4\pmod 2.$$

Computations suggest:

\begin{conj}\label{mod3}
We have
$$c_k(t)\equiv (t+1)^{3k}\det (tI+P(k))\pmod 3\ .$$
\end{conj}

This conjecture, together with Theorem  \ref{ppowerformula} yields
conjectural recursive formulas for
$p_n(t)=\det (tI(n)-P(n))\pmod 3$ as follows: Set
$p_0(t)=1\pmod 3,\ p_1(t)=1-t\pmod 3$. For $n=3^l\pm k>1$ 
with $0\leq k<\frac{3^l}{2}$
the characteristic polynomial $\chi_n(t)\pmod 3$ is then conjecturally
given by
$$\begin{array}{ll}
\displaystyle (t-1)^{3^{l}-3k}\ \det(t^2I+P(k))\qquad&\displaystyle
\hbox{if }n=3^l-k\ ,\cr
\displaystyle (t-1)^{3^{l}-3k}\ 
(t+1)^{3k}\ \det(tI+P(k))\qquad&\displaystyle
\hbox{if }n=3^l+k\ .\end{array}$$
In particular, all roots of $\chi_n(t)$ modulo~$3$ should be of
multiplicative order a power of $2$ in the algebraic closure of $\F_3$.

%

We conclude finally by mentioning a last conjectural observation:

\begin{conj} Given a prime-power $q=p^l\equiv 2\pmod 3$, we have
$$\chi_{(q+1)/3}(t)\equiv (t+1)^{(q+1)/3}\pmod p$$
and 
$$\chi_{(2q-1)/3}(t)\equiv(t+1)^{(q+1)/3}\ (t-1)^{(q-2)/3}\pmod p.$$
\end{conj}

\begin{rem} (i)
The matrix 
$C=P(\frac{q+1}{3})+I(\frac{q+1}{3})$ for $q=p^l\equiv2\pmod 3$ a 
prime-power, appears to have a unique Jordan block
of maximal length over~$\F_p$. If so,
the rows of $C^{(q+1)/6}$ generate
a self-dual code over~$\F_p$.

\ \ (ii) Given a prime power $q=p^l\equiv 2\pmod 3$ as above we set
$n=\frac{2q+2}{3}$ and $k=\frac{2q-1}{3}$. 
We conjecture that the characteristic polynomial
of the matrix $\tilde P_k(n)$ with coefficients
$$\tilde p_{i,j}={i+j+2k\choose i+k},\ 0\leq i,j<n$$
satisfies
$\det(tI-\tilde P_k(n))\equiv (1+t)^n\pmod p$.
\end{rem}

\begin{rem}
In \cite[Theorems 32 and 35]{K} Krattenthaler gives evaluations
of determinants related to ours, namely of $\det(\om I+Q(n))$
where $\om$ is a sixth root of unity, and $Q(n)$ has entries
${2\mu+i+j\choose j}$ $(0\le i,j<n)$.
\end{rem}

The sequel of this paper is organized as follows:

Section 2 is devoted to autosimilar matrices. Such matrices
generalize the matrix ${\overline P}(\infty)_2$ and their properties
imply easily Theorem \ref{detmod2}.

Section 3 contains proofs of Proposition \ref{order3forprimep} 
and Theorem \ref{ppowerformula}.

Section 4 contains proofs of Theorems~\ref{mod2} and
\ref{gammaformula}.

\section{Autosimilar matrices}

Let $b\geq 1$ be a natural integer. An infinite matrix
$M$ with coefficients $m_{i,j}$ ($i,j\ge0$) is
{\it $b$-autosimilar} if $m_{0,0}=1$ and if 
$$m_{s,t}=\prod_i m_{\sigma_i,\tau_i}$$
where the indices $s=\sum \sigma_i b^i,\ t=\sum \tau_i b^i$ are written in base
$b$, that is, $\sigma_i,\tau_i\in\{0,\dots,b-1\}$ for all $i=0,1,2,\dots$.

We denote by $M(n)$ the finite sub-matrix of $M$ with coefficients
$m_{i,j},\ 0\leq i,j<n$. A $b$-autosimilar matrix $M$ is
{\it non-degenerate} if the determinants
$$\det (M(n))$$ are invertible for $n=2,\dots,b$.  

\begin{thm}\label{LDU} Let $b\geq 2$ be an integer and let $M$ be a 
$b$-autosimilar matrix which is non-degenerate. One has then a factorization
$$M=LDU$$
where $L,D,U$ are $b$-autosimilar and where
$L$ is unipotent lower-triangular, $D$ is diagonal and $U$ is unipotent
upper-triangular.
\end{thm}

\begin{cor}\label{detautosimilar} 
Given a non-degenerate $b$-autosimilar
matrix $M$ one has
$$\det (M(n))=\prod_{i=0}^{n-1} d_{\nu_i}$$
for all $n=\sum \nu_i b^i$ with $d_0=1$ and
$$d_k=\det (M(k+1))/\det (M(k))$$
for $k=1,\dots,b-1$.
\end{cor}

\begin{rem}\label{remautosimilar} 

In general, one can compute determinants of arbitrary
$b$-autosimilar matrices over a field $K$
by applying Corollary \ref{detautosimilar} to 
the $b$-autosimilar matrix obtained from a generic perturbation of the
form $$M_t(b)=(1-t)M(b)+tP(b)$$
(where $P(b)$ is a suitable matrix) and working over the rational function
field~$K(t)$.
\end{rem}

{\bf Proof of Theorem \ref{LDU}.} The genericity of $M$ implies that
$$M(b)=L(b)D(b)U(b)$$
where $L(b)$ and $U(b)$ are unipotent upper and lower triangular matrices
and the diagonal matrix $D(b)$ has entries $d_{0,0}=1$
and $d_{k,k}=\det (M(k+1))/\det (M(k))$
for $k=1,\dots,b-1$. Extending $L(b)$, $D(b)$ and $U(b)$ in the unique
possible way to infinite $b$-autosimilar matrices $L$, $D$ and $U$ we have
$$\begin{array}{ll}
\displaystyle (LDU)_{s,t}
&\displaystyle =\sum_k L_{s,k}D_{k,k}U_{k,t}\cr
&\displaystyle =\sum_{k=\sum\kappa_i b^i} \prod_i
L_{\sigma_i,\kappa_i}D_{\kappa_i,\kappa_i}U_{\kappa_i,\tau_i}\cr
&\displaystyle =\prod_i\sum_{\kappa_i=0}^{b-1}
L_{\sigma_i,\kappa_i}D_{\kappa_i,\kappa_i}U_{\kappa_i,\tau_i}\cr
&\displaystyle =\prod_i M_{\sigma_i,\tau_i}=M_{s,t}\end{array}$$
for all $s=\sum \sigma_ib^i,t=\sum \tau_ib^i\in{\bf N}$. \hfill $\Box$

The identity
$$\det (M(n))=\det (D(n))$$
implies immediately Corollary \ref{detautosimilar}. 

\subsection{Binomial coefficients modulo a prime $p$}

Let $p$ be a prime number. We have then
$$(1+x)^n=\prod (1+x)^{\nu_i p^i}\equiv (1+x^{p^i})^{\nu_i}\pmod p$$
(using properties of the Frobenius automorphism in characteristic $p$).
This implies 
immediately the equality
$${n\choose k}=\prod_i {\nu_i\choose \kappa_i}$$
allowing (for small primes)
an efficient computation of binomial coefficients$\pmod p$.

This equality shows that the reductions modulo 2 or 3 of the
symmetric Pascal triangle $P$ with coefficients 
$${\overline p}_{i,j}=\left({i+j\choose i}\pmod 2\right)\in\{0,1\}$$
respectively
$${\overline p}_{i,j}=\left({i+j\choose i}\pmod 3\right)\in\{-1,0,1\}$$
are $2-$ (respectively $3-$) autosimilar matrices.

For $p=2$ we have
$$\left(\begin{array}{cc}
1&1\cr 1&0\end{array}\right)=
\left(\begin{array}{cc}
1&0\cr 1&1\end{array}\right)
\left(\begin{array}{rr}
1&0\cr 0&-1\end{array}\right)
\left(\begin{array}{rr}
1&1\cr 0&1\end{array}\right)$$
which yields $d_0=1,d_1=-1$ and Corollary \ref{detautosimilar} 
implies now Theorem \ref{detmod2}.

\begin{rem}\label{inversemod2} 
One can show that the inverse of the integral matrix
${\overline P}(n)_2$ considered in Theorem \ref{detmod2} has all its
coefficients in $\{-1,0,1\}$ for all $n$.
\end{rem}

For $p=3$ we have
$$\left(\begin{array}{rrr}
1&1&1\cr
1&-1&0\cr
1&0&0\end{array}\right)=
\left(\begin{array}{rrr}
1&0&0\cr
1&1&0\cr
1&\frac{1}{2}&1\end{array}\right)
\left(\begin{array}{rrr}
1&0&0\cr
0&-2&0\cr
0&0&-\frac{1}{2}\end{array}\right)
\left(\begin{array}{rrr}
1&1&1\cr
0&1&\frac{1}{2}\cr
0&0&1\end{array}\right)$$
This shows that $\det (\overline{P}(n)_3)$ (over $\bf Z$) equals
$(-2)^{a-b}$ where $a$ and $b$ are the number of digits $1$ and $2$ needed
in order to write all natural integers $<n$ in base $3$.

\section{Proofs of Proposition \ref{order3forprimep}
and Theorem \ref{ppowerformula}}

\begin{pf}{\bf of Proposition \ref{order3forprimep}}
Let $R$ be a commutative ring, and let
$$A=\mx{a}{b}{c}{d}\in\GL(2,R).$$
Then $A$ determines a (graded $R$-algebra) automorphism $\phi_A$ of $
R[X,Y]$
via $\phi_A(X)=aX+bY$ and $\phi_A(Y)=cX+dY$, or alternatively
$$\vc{\phi_A(X)}{\phi_A(Y)}=A\vc{X}{Y}.$$
It is easy to see that $\phi_A\circ\phi_B=\phi_{BA}$.
Each $\phi_A$ restricts to an $R$-module automorphism of the homogeneous
polynomials $R[X,Y]_{n-1}$ of degree $n-1$.
Let $A^{(n)}$ denote the matrix of this endomorphism with respect to the
basis $X^{n-1}$, $X^{n-2}Y$, $X^{n-3}Y^2,\ldots,Y^{n-1}$, that is
$$\left(
\begin{array}{c}
{\phi_A(X^{n-1})}\\
{\phi_A(X^{n-2}Y)}\\
{\phi_A(X^{n-3}Y^2)}\\
\vdots\\
{\phi_A(Y^{n-1})}\\
\end{array}
\right)
=A^{(n)}
\left(\begin{array}{c}
{X^{n-1}}\\
{X^{n-2}Y}\\
{X^{n-3}Y^2}\\
\vdots\\
{Y^{n-1}}\\
\end{array}
\right).
$$
Then $A^{(n)}\in\GL(n,R)$ and $(AB)^{(n)}=A^{(n)}B^{(n)}$.
(Another way of expressing this is to say that $A^{(n)}$ is the
$(n-1)$-th symmetric power of~$A$.)

Let us specialize to the case $R=\F_p=\Z/p\Z$ and $n=p^l$.
In this case $A^{(n)}=I$ if and only if $A$ is a scalar matrix.
The matrix
$$A=\mx{1}{-1}{1}{0}$$
yields $A^{(n)}\equiv P(p^{l})\pmod p$. Since $A^3=-I$, the matrix
$A^{(n)}$ has order~$3$.

Let us now compute the multiplicities of the three eigenvalues of
$P=P(p)\pmod p$ over $\F_p$ (the formula for $P(p^l)$ is
then a straightforward consequence of the fact the $P(p^l)$ is
the $l-$fold Kronecker product of $P(p)$ with itself).

The easy identity ${2k\choose k}={(p-1)/2\choose k}(-4)^k\pmod p$
for $p$ an odd prime and $0\leq k\leq (p-1)/2$ shows
$$\sum_{k=0}^{(p-1)/2}{2k\choose k}\left(\frac{-x}{4}\right)^k\equiv
(1+x)^{(p-1)/2}\pmod p$$
and yields $\hbox{tr}(P)\equiv(-3)^{(p-1)/2}\equiv 
\epsilon(p)\pmod p$ (where $\epsilon(p)\in
\{-1,0,1\}$ satisfies $\epsilon(p)\equiv p\pmod 3$) by quadratic reciprocity.

Since the characteristic polynomial for $P$ has antisymmetric
coefficients ($\alpha_k=-\alpha_{p-k}$) the two eigenvalues
$\not=1$ of $P$ have equal multiplicity $r$.
Lifting into positive integers $\leq \frac{p-1}{2}$ the solution
of the linear system $-r+(p-2r)\equiv\hbox{tr}(P)\pmod p$
yields now the result.

The case $p=2$ is easily solved by direct inspection.
\end{pf}

\begin{rem}\label{groups}
Recall that we have (with the notations of the above proof) 
$P=P(n)=A^{(n)}\pmod p$ for $n=p^l$ and introduce $L=L(n)=B^{(n)}\pmod p$
and $\tilde L=\tilde L(n)=C^{(n)}\pmod p$ where
$$A=\mx{1}{-1}{1}{0},B=\mx{1}{0}{-1}{-1},C=\mx{1}{0}{1}{-1}.$$

It is straightforward to check that $L$ and $\tilde L$ have coefficients
$$l_{i,j}=(-1)^i{i\choose j}\pmod p\qquad \hbox{ and }\qquad
\tilde l_{i,j}=(-1)^j{i\choose j}\pmod p$$
for $0\leq i,j<n$.

Then $A^3=-I$, but $(-I)^{(n)}$
is the identity. Hence $P^3=I$. Also $C^2=I$ and $CAC=A^{-1}$.
It follows that $A$ and $C$ generate a dihedral group of order 12,
containing $-I$. Hence $A^{(n)}=P$ and $C^{(n)}=\ti{L}$ generate
a dihedral group of order~6. 

The group $G_p$ generated by $P$ and $L$ depends on the prime $p$ (but not
on the power $l$ of $n=p^l$). It is isomorphic to a subgroup of
$\hbox{PGL}_2(\F_p)$. For all but finitely many primes~$p$, $G_p$ is
isomorphic to $\hbox{PSL}_2(\F_p)$ or $\hbox{PGL}_2(\F_p)$ according
to whether $-1$ is or is not a square in $\F_p$.
The exceptional primes are $5$, $7$ and $29$
where $G_p$ has order $24$, $42$ and $120$ respectively.
\end{rem}

\begin{pf}{\bf of Theorem \ref{ppowerformula}}
Using Proposition \ref{order3forprimep},
we can rewrite the equation to be proved as
$$(t^3-1)^k \det(tI-P(q-k))
\equiv\det(tI-P(q))\det(t^2I+P(k))\pmod p.$$
Here, and in the sequel, we write $I$ for $I(n)$ whenever this
notation is unambiguous; also we denote the zero matrix of any size by~$O$.

We now work over the field~$\F_p$. Unless otherwise stated
vectors will be row vectors.

It is convenient to define a category $\ce=\ce_{\F_p}$
as follows. Its objects will be pairs $(V,\al)$
where $V$ is a finite-dimensional vector space over $\F_p$
and $\al$ is a vector space endomorphism of~$V$. A morphism
$\phi:(V,\al)\to(W,\be)$ in $\ce$ will be a linear map $\phi:V\to W$
with $\phi\circ\al=\be\circ\phi$. (In fact $\ce$ is equivalent to
the category of finitely generated torsion modules over the
polynomial ring $\F_p[X]$.)
If $(V,\al)$ is an object of $\ce$
we define $\chi(V,\al,t)$ as the characteristic polynomial of $\al$
acting on~$V$, that is, $\chi(V,\al,t)=\det(tI-A)$ where $A$ is a matrix
representing $\al$ with respect to some basis of~$V$.
An $r$ by $r$ matrix $A$ defines an object $((\F_p)^r,\al)$, denoted by
$((\F_p)^r,A)$, where $\al$ is the endomorphism defined by~$A$.

It is easy to see that $\ce$ is an abelian category, and that if
$$0\to(V,\al)\to(X,\ga)\to(W,\be)\to0$$
is a short exact sequence, then $\chi(X,\ga,t)=\chi(V,\al,t)\chi(W,\be,t)$.
This is because there is a basis for $X$ with respect to which the
matrix of $\ga$ (acting on row vectors from the the right) is
$$\left(\begin{array}{cc}
A&O\\
C&B
\end{array}
\right)$$
where $A$ and $B$ are matrices representing $\al$ and $\be$ respectively.

Set $k'=q-k$. We can partition the Pascal matrices $P(k')$
and $P(q)$ as follows:
$$P(k')=\left(\begin{array}{cc}
A&B\\
B^t&C
\end{array}
\right)\qquad\textrm{and}\qquad
P(q)=\left(\begin{array}{ccc}
A&B&D\\
B^t&C&O\\
D^t&O&O
\end{array}
\right)$$
where $A=P(k)$.

Let $\ov{A}$ denote the matrix obtained by rotating $A$ through $180^\circ$.
Then $P(q)^2=\ov{P(q)}$ and $P(q)^3=I$. Hence
$$P(q)^2=\left(\begin{array}{ccc}
O&O&\ov{D^t}\\
O&\ov{C}&\ov{B^t}\\
\ov{D}&\ov{B}&\ov{A}
\end{array}
\right).$$
Thus
$$A^2+BB^t+DD^t=O$$
and so
$$P(k')^2=
\left(\begin{array}{cc}
-DD^t&O\\
O&\ov{C}
\end{array}
\right).$$
 
 From $P(q)^2=\ov{P(q)}$ it follows that $AD=\ov{D^t}$ and from
$\ov{P(q)}P(q)=I$ it follows that $\ov{D^t}D^t=I$. Hence $ADD^t=I$
and so
$$P(k')^2=
\left(\begin{array}{cc}
-A^{-1}&O\\
O&\ov{C}
\end{array}
\right).$$

Let $V=(\F_p)^q$ and $X=(\F_p)^{3k}$. Let
$$Q_1=\left(\begin{array}{ccc}
O&I(k)&O\\
O&O&I(k)\\
I(k)&O&O
\end{array}
\right).$$
Let $\phi:X\to V$ be the map defined by the matrix
$$\left(\begin{array}{ccc}
I&O&O\\
A&B&D\\
O&O&\ov{D^t}
\end{array}
\right).$$
Then
$$Q_1
\left(\begin{array}{ccc}
I&O&O\\
A&B&D\\
O&O&\ov{D^t}
\end{array}
\right)
=\left(\begin{array}{ccc}
A&B&D\\
O&O&\ov{D^t}\\
I&O&O
\end{array}\right)
$$
and
$$\left(\begin{array}{ccc}
I&O&O\\
A&B&D\\
O&O&\ov{D^t}
\end{array}
\right)
P(q)=
\left(\begin{array}{ccc}
I&O&O\\
A&B&D\\
O&O&\ov{D^t}
\end{array}
\right)
\left(\begin{array}{ccc}
A&B&D\\
B^t&C&O\\
D^t&O&O
\end{array}\right)
=\left(\begin{array}{ccc}
A&B&D\\
O&O&\ov{D^t}\\
I&O&O
\end{array}
\right)$$
where we have used the formulas $P(q)^2=\ov{P(q)}$ and $\ov{P(q)}P(q)=I$.
Hence $\phi$ is a morphism from $((\F_p)^{3k},Q_1)$
to $((\F_p)^q,P(q))$ in~$\ce$.

Let $W=(\F_p)^{k'}$ and $Y=(\F_p)^{2k}$. Let
$$Q_2=\left(\begin{array}{cc}
O&I(k)\\
-A^{-1}&O\\
\end{array}
\right).$$
Let $\psi:Y\to W$ be the map defined by the matrix
$$\left(\begin{array}{cc}
I&O\\
A&B
\end{array}
\right).$$
Then
$$Q_2
\left(\begin{array}{cc}
I&O\\
A&B
\end{array}
\right)
=\left(\begin{array}{cc}
A&B\\
-A^{-1}&O
\end{array}\right)
$$
and
$$\left(\begin{array}{cc}
I&O\\
A&B
\end{array}
\right)P(k')=
\left(\begin{array}{cc}
I&O\\
A&B
\end{array}
\right)
\left(\begin{array}{cc}
A&B\\
B^t&C
\end{array}
\right)
=\left(\begin{array}{cc}
A&B\\
-A^{-1}&O
\end{array}\right)
$$
where we have used the formula
$$P(k')^2=
\left(\begin{array}{cc}
-A^{-1}&O\\
O&\ov{C}
\end{array}\right).
$$
Hence $\psi$ is a morphism from $((\F_p)^{2k},Q_2)$
to $((\F_p)^{k'},P(k'))$ in~$\ce$.

We need to divide into the cases $k\le q/3$ and $k\ge q/3$.
In the former cases $\phi$ and $\psi$ are injective
and in the latter case they are surjective. In the former
case we consider their cokernels, in the latter case their kernels.

The matrix $B$ has size $k$ by~$q-2k$. If $B$ has rank $k$
(which is only possible if $k\le q/3$) then $\phi$ and $\psi$ are
injective. If $B$ has rank $q-2k$
(which is only possible if $k\ge q/3$) then $\phi$ and $\psi$ are
surjective.

The matrix $B$ contains a submatrix
$$\left({i+j+k\choose i}\right)_{i,j=0}^{r-1}$$
where $r=\min(k,q-2k)$. This submatrix has
determinant~$1$ (consider it as a matrix over $\Z$
and reduce it to a Vandermonde matrix or see for instance \cite{B}). 
Thus $B$ has rank $r$
and indeed $\phi$ and $\psi$ are injective for $k\le q/3$
and surjective for $k\ge q/3$.

Consider first the case where $k\le q/3$. Let $(X_1,\tht_1)$
and $(X_2,\tht_2)$ denote the cokernels of 
$\phi:((\F_p)^{3k},Q_1)\to((\F_p)^q,P(q))$
and $\psi:((\F_p)^{2k},Q_2)\to((\F_p)^{k'},P(k'))$ in $\ce$.
Then
$$\chi((\F_p)^q,P(q),t)=\chi((\F_p)^{3k},Q_1,t)\chi(X_1,\tht_1,t)$$
and
$$\chi((\F_p)^{k'},P(k'),t)=\chi((\F_p)^{2k},Q_2,t)\chi(X_2,\tht_2,t).$$
It is apparent that
$$\chi((\F_p)^{3k},Q_1,t)=(t^3-1)^k$$
and
$$\chi((\F_p)^{2k},Q_2,t)=\det(t^2 I+A^{-1})=\det(t^2 I+A)$$
as $A$ and $A^{-1}$ are similar. Hence
$$\det(tI-P(q))=(t^3-1)^k\chi(X_1,\tht_1,t)$$
and
$$\det(tI-P(k'))=\det(t^2I+A)\chi(X_2,\tht_2,t).$$
It suffices to prove that $(X_1,\tht_1)$ and
$(X_2,\tht_2)$ are isomorphic in~$\ce$.

As $\ov{D^t}$ is nonsingular, it is apparent that $X_1$
is isomorphic to $(\F_p)^{q-2k}/Y$ where $Y$ is the row space of
$B$ and that the action of $\tht_1$ is induced by that of the matrix~$C$
on $(\F_p)^{q-2k}$. It is even more apparent that $X_2$
is isomorphic to $(\F_p)^{q-2k}/Y$ and
that the action of $\tht_2$ is induced by~$C$. Hence
$(X_1,\tht_1)$ and $(X_2,\tht_2)$ are isomorphic in~$\ce$.
This completes the argument in the case $k\le q/3$.

Now suppose that $k\ge q/3$.
Let $(K_1,\tht_1)$ and $(K_2,\tht_2)$ denote the kernels of 
$\phi:((\F_p)^{3k},Q_1)\to((\F_p)^q,P(q))$
and $\psi:((\F_p)^{2k},Q_2)\to((\F_p)^{k'},P(k'))$ in $\ce$.
Then
$$\chi((\F_p)^q,P(q),t)\chi(K_1,\tht_1,t)=\chi((\F_p)^{3k},Q_1,t)$$
and
$$\chi((\F_p)^{k'},P(k'),t)\chi(K_2,\tht_2,t)=\chi((\F_p)^{2k},Q_2,t).$$
Hence
$$\frac{(t^3-1)^k}{\det(tI-P(q))}=\chi(K_1,\tht_1,t)$$
and
$$\frac{\det(t^2I+A)}{\det(tI-P(k'))}=\chi(K_2,\tht_2,t).$$
It suffices to prove that $(K_1,\tht_1)$ and
$(K_2,\tht_2)$ are isomorphic in~$\ce$.

As $\ov{D^t}$ is nonsingular and has inverse $D^t$, it is apparent that
$$K_1=\{(-uA,u,-uDD^t)=(-uA,u,-uA^{-1}):u\in(\F_p)^k,uB=0\}$$
and we have
$$(-uA,u,-uA^{-1})Q_1=(-uA^{-1},-uA,u)\ .$$
Also
$$K_2=\{(-uA,u):u\in(\F_p)^k,uB=0\}$$
and
$$(-uA,u)Q_2=(-uA^{-1},-uA)\ .$$
Hence the linear map
$$(-uA,u,-uA^{-1})\longmapsto(-uA,u)$$
induces an isomorphism between $(K_1,\tht_1)$ and $(K_2,\tht_2)$.
\end{pf}

\section{Proofs for the prime $p=2$}

\begin{pf}{\bf of Theorem \ref{mod2}.} Set $n=2^l-k$ and $q=2^l$ where 
$1\leq k\leq 2^{l-1}$.

Theorem \ref{ppowerformula} yields then over $\F_2$
$$\chi_n(t)= \chi_{q-k}(t)=
(t^2+t+1)^{(q-\epsilon(q))/3-k}(t+1)^{(q+2\epsilon(q))/3-k}
\det(tI+P(k))^2$$
since $x\longmapsto x^2$ is an automorphism in characteristic $2$.

By induction on $l$, the only possible irreducible factors of 
$\det (tI(n)-P(n))\pmod 2$ are $(1+t)$ and $(1+t+t^2)$. The 
multiplicity $\mu(n)=\mu(2^l-k)$ of the factor $(1+t)$ in this
polynomial is hence recursively defined by
$$\mu(n)=\frac{2^l+2(-1)^l}{3}-k+2\mu(k)$$
and coincides hence with the sequence $\gamma$ of Theorem  \ref{mod2}.
The remaining factor of $\det (tI(n)-P(n))\pmod 2$ is hence given by
$(1+t+t^2)^{\gamma_2(n)}$ where $\gamma_2(n)=\frac12(n-\gamma(n))$
and this proves the result.
\end{pf}

\begin{pf}{\bf of Theorem  \ref{gammaformula}.} 
%
We have for $0\leq k\leq 2^{l-1}$
\begin{eqnarray*}
\gamma(2^l+k)&=&\gamma(2^{l+1}-(2^l-k))\\
&=&\frac{2^{l+1}-2(-1)^l}{3}-2^l+k+2\gamma(2^l-k)\\
&=&\frac{2^{l+1}-2(-1)^l}{3}-2^l+k+2\frac{2^{l}+2(-1)^l}{3}-2k+4\gamma(k)
\end{eqnarray*}
which is assertion (i).

We have for all $2^{l-2}\leq k\leq 2^{l-1}$
\begin{eqnarray*}
\gamma(2^l-k)&=&\frac{2^{l}+2(-1)^{l}}{3}-k+\gamma(k)+
\gamma(2^{l-1}-(2^{l-1}-k))\\
&=&\frac{2^{l}+2(-1)^{l}}{3}-k+\gamma(k)+
\frac{2^{l-1}-2(-1)^l}{3}-2^{l-1}+k+2\gamma(2^{l-1}-k)\\
&=&\gamma(k)+2\gamma(2^{l-1}-k)
\end{eqnarray*}
which proves assertion (ii).

Similarly, we have for $1\leq k\leq 2^l$
\begin{eqnarray*}
\gamma(2^l+k)-\gamma(2^l+k-1)&=&
\gamma(2^{l+1}-(2^l-k))-\gamma(2^{l+1}-(2^l-k+1))\\
&=&1+2\gamma(2^l-k)-2\gamma(2^l-k+1)
\end{eqnarray*}
which proves assertion (iii).

Writing $2n=2^l-2k$ with $1\leq k\leq2^{l-2}$ we have,
using induction on $n$,
\begin{eqnarray*}
\gamma(2^l-2k)&=&\frac{2^l-(-1)^l}{3}-2k+2\gamma(2k)\\
&=&\frac{2^l-(-1)^l}{3}-2k+2\left(k-\gamma(k)\right)\\
&=&\left(2^{l-1}-k\right)-\left(\frac{2^{l-1}-(-1)^{l-1}}{3}-k+2\gamma(k)
\right)\\
&=&\left(2^{l-1}-k\right)-\gamma(2^{l-1}-k)
\end{eqnarray*}
which proves the first equality of assertion (iv) (this equality follows 
also from the fact that $P(2n)$ is the Kronecker product of $P(n)$ with $P(2)$
over~$\F_2$).

The second identity of assertion (iv) amounts to the equality
$$\gamma(2n-1)-\gamma(2n)=\frac{4^{b(2n-1)}-1}{3}\ .$$
We prove first by induction on $n$ that this identity is
equivalent to the last identity.

The last identity and induction yield
\begin{eqnarray*}
\gamma(2n-1)-\gamma(2n)&=&\gamma(2n-1)-
\gamma(2n-2)+\gamma(2n-2)-\gamma(2n)\\
&=&\frac{2^{1+2b(n-1)}+1}{3}-1+\gamma(n)-\gamma(n-1).
\end{eqnarray*}
We now divide into cases according to the parity of~$n$.

Suppose first that $n=2m$ is even. Then inductively
$$\gamma(n)-\gamma(n-1)=\gamma(2m)-\gamma(2m-1)=-\frac{4^{b(2m-1)-1}}{3}
=-\frac{4^{b(n-1)-1}}{3}$$
Hence
$$\gamma(2n-1)-\gamma(2n)=-1+\frac{2^{1+2b(n-1)}+1}{3}
-\frac{2^{2b(n-1)}-1}{3}=\frac{2^{2b(n-1)}-1}{3}.$$
But
$$2^{2b(n-1)}=4^{b(n-1)}=4^{b(2n-1)}$$
as the binary representation of $n-1$ ends in $1$ and that of $2n-1$
is obtained by appending~$1$.

Now suppose that $n=2m+1$ is odd. Then
$$\gamma(n)-\gamma(n-1)=\gamma(2m+1)-\gamma(2m)=\frac{2^{1+2 b(m)}+1}{3}
=\frac{2^{1+2 b(2m)}+1}{3}.$$
Hence
$$\gamma(2n-1)-\gamma(2n)=-1+\frac{2^{1+2b(n-1)}+1}{3}
+\frac{2^{1+2b(n-1)}+1}{3}=\frac{2^{2+2b(n-1)}-1}{3}.$$
But
$$2^{2+2b(n-1)}=4^{1+b(n-1)}=4^{b(2n-1)}$$
as the binary representation of $n-1$ ends in $0$ and that of $2n-1$
is obtained by appending~$1$.

This completes the proof of equivalence of the two last identities 
in assertion (iv).

We prove now the last identity by induction on $n$.

The last identity of assertion (iv) is equivalent to
$$\gamma(2n+1)-\gamma(2n)=\frac{2^{1+2b(n)}+1}{3}\ .$$
Writing $2n+1=2^l+k$ with $1\leq k<2^l$ and applying assertion (iii)
and the second identity of assertion (iv) (which holds by induction)
we have
\begin{eqnarray*}
\gamma(2n+1)-\gamma(2n)&=&1+2\gamma(2^l-k)-2\gamma(2^l+1-k)\\
&=&1+2\frac{4^{b(2^l-k)}-1}{3}\\
&=&\frac{2^{1+2b(2^l-k)}+1}{3}
\end{eqnarray*}
Since $(2^l+k-1)+(2^l-k)=2^{l+1}-1$ and since $2^l+k-1$ is even and greater
than $2^l-k$, they have the same number of blocks $1\dots1$ in their
binary expansion. This shows $b(2^l-k)=b(2n)=b(n)$ and establishes the
last identity of assertion (iv). 
\end{pf}

The first author wishes to thank J.-P.~Allouche, F.~Sigrist, U.~Vishne
and  A.~Wassermann for interesting comments and remarks.

Roland Bacher, Institut Fourier, UMR 5582,
Laboratoire de Math\'ematiques, BP 74, 38402 St. Martin d'H\`eres Cedex,
France, Roland.Bacher@ujf-grenoble.fr

Robin Chapman, University of Exeter, School of Mathematical Sciences,
North Park Road, EX4 4QE Exeter, UK, rjc@maths.ex.ac.uk

\end{document}